\documentclass[12pt]{amsart}

\usepackage{amsmath,amsthm,amscd,euscript}
\def \r{\mathbb R}

\def \z{\mathbb Z}

\theoremstyle{remark}

\theoremstyle{definition}

\newtheorem{example}{Example}

\title[Three-dimensional continued fractions]{Three examples of three-dimensional continued fractions
in the sense of Klein.}
\author{Oleg~Karpenkov}
\date{19 December 2005}
\thanks{Partially supported
by RFBR grant SS-1972.2003.1 and by RFBR grant 05-01-01012a.}

\keywords{Integer lattice, continued fractions, convex hulls}

\email[Oleg Karpenkov]{karpenk@ceremade.dauphine.fr}

\address{CEREMADE - UMR 7534
Universit\'e Paris-Dauphine France -- 75775 Paris SEDEX 16}

\begin{document}
\input epsf
\maketitle

The problem of investigation of the simplest $n$-dimensional
continued fraction for $n\ge 2$ was posed by V.~Arnold. The answer
for the case of $n=2$ can be found in the works of
E.~Korkina~\cite{Kor2} and G.~Lachaud~\cite{Lac2}. In present
work we study the case of $n=3$. The author is grateful to
V.~Arnold for constant attention to this work and useful remarks
and Universit\'e Paris-Dauphine
--- CEREMADE  for the hospitality and excellent working
conditions.

{\bf Definitions.} A point of $\r^{n{+}1}$ is called {\it integer}
if  all its coordinates are integers. A hyperplane is called {\it
integer} if all its integer vectors generate $n$-dimensional
sublattice of integer lattice. Consider some integer hyperplane
and an integer point in the complement to this plane. Let the
Euclidean distance from the given point to the given plane equals
$l$.  The minimal value of nonzero Euclidean distances from
integer points of the space $\r^{n+1}$ to the plane is denoted by
$l_0$. The ratio $l/l_0$ is said to be the {\it integer distance}
from the given integer point to the given integer hyperplane.

{\bf Definition of multidimensional continued fraction in the
sense of Klein.} Consider arbitrary $n{+}1$ hyperplanes
in~$\r^{n+1}$ that intersect at the unique point: at the origin.
Assume also that all the given planes do not contain any integer
point different to the origin. The complement to these hyperplanes
consists of $2^{n+1}$ open orthants. Consider one of the orthants.
The boundary of the convex hull of all integer points except the
origin in the closure of the orthant is called the {\it sail} of
the orthant. The set of all $2^{n+1}$ sails is called the {\it
$n$-dimensional continued fraction} constructed accordingly to the
given $n{+}1$ hyperplanes. Two $n$-dimensional continued
fractions are said to be {\it equivalent} if there exist a linear
lattice preserving transformation of $\r^{n+1}$ taking all sails
of one continued fraction to the sails of the other continued
faction.

We associate to any hyperbolic irreducible operator $A$ of
$SL(n{+}1,\z)$ an $n$-dimensional continued fraction constructed
accordingly to the set of all $n{+}1$ eigen-hy\-per\-pla\-nes for
$A$. Any sail of such continued fraction is homeomorphic to
$\r^n$. From Dirichlet unity theorem it follows that the group of
all $SL(n{+}1,\z)$-operators commuting with $A$ and preserving
the sails is homeomorphic to $\z^n$ and its action is free (we
denote this group by $\Xi(A)$). A {\it fundamental domain} of the
sail with respect to the action of the group $\Xi(A)$ is a face
union that contains exactly one face of the sail from each orbit.
(For more information see~\cite{Arn2},~\cite{Kar1},~\cite{Kor2}
and~\cite{Mou2}.)

Denote by $A_{a,b,c,d}$ the following integer operator
% Sylvester matrix with charpoly ...
$$
\left(
\begin{array}{cccc}
0&1&0&0\\
0&0&1&0\\
0&0&0&1\\
a&b&c&d\\
\end{array}
\right)
.
$$

\begin{example} Consider the operator $A_1=A_{1,-3,0,4}$. The group
$\Xi(A_1)$ is generated by the operators $B_{11}=A_1^{-2}$,
$B_{12}=(A_1{-}E)^2A_1^{-2}$, and
$B_{13}=(A_1{-}E)^2(A_1{+}E)A_1^{-2}$. Let us enumerate all
three-dimensional faces for one of the fundamental domains of the
sail containing the vertex $(0,0,0,1)$. Let $V_{10}=(-3, -2, -1,
1)$, $V_{1,4i+2j+k}=B_{11}^i B_{12}^j B_{13}^k (V_{10})$ for
$i,j,k \in \{0,1\}$. One of the fundamental domains of the sail
contains the following three-dimensional faces:
$T_{11}=V_{10}V_{12}V_{14}V_{15}$,
$T_{12}=V_{12}V_{14}V_{15}V_{16}$,
$T_{13}=V_{12}V_{15}V_{16}V_{17}$,
$T_{14}=V_{12}V_{13}V_{15}V_{17}$,
$T_{15}=V_{10}V_{12}V_{13}V_{15}$,
$T_{16}=V_{10}V_{11}V_{13}V_{15}$, and
$T_{17}=V_{10}V_{11}V_{12}V_{13}$ (see on the figure from the
left). All listed tetrahedra are taken by some integer affine
transformations to the unit basis tetrahedron. The integer
distance from the origin to the planes containing the faces
$T_{11},\ldots, T_{17}$ equal 4, 3, 2, 4, 3, 2, and 1,
respectively.
\end{example}

{\bf Statement 1.} The continued fraction constructed for any
hyperbolic matrix of $SL(4,\z)$ with irreducible characteristic
polynomial over rationals and with the sum of absolute values of
the elements smaller than 8 is equivalent to the continued
fraction of Example~1.

{\bf Statement 2.} The symmetry (not commuting with $A_1$)
defined by the matrix
$$
\left(
\begin{array}{rrrr}
4&-16&17&-3\\
3&-11&11&-2\\
3&-8&6&-1\\
6&-8&-2&1\\
\end{array}
\right)
$$
acts on the sail of Example~1. This symmetry permutes the
equivalence classes (with respect to the action of $\Xi(A_1)$) of
tetrahedra $T_{11}$ and $T_{14}$, $T_{12}$ and $T_{15}$, $T_{13}$
and $T_{16}$, and takes the class of $T_{17}$ to itself.

\begin{example} Let us consider the operator $A_2=A_{1,-4,1,4}$. The group
$\Xi(A_2)$ is generated by the operators $B_{21}=A_2^{-2}$,
$B_{22}=(A_2{-}E)^2A_2^{-2}$, and $B_{23}=(A_2{+}E)A_2^{-1}$. Let
us enumerate all three-dimensional faces for one of the
fundamental domains of the sail containing the vertex
$(0,0,0,1)$. Let $V_{20}=(-4, -3, -2, 0)$,
$V_{2,4i+2j+k}=B_{21}^i B_{22}^j B_{23}^k (V_{20})$ for $i,j,k
\in \{0,1\}$. One of the fundamental domains of the sail contains
the following three-dimensional faces:
$T_{21}=V_{20}V_{21}V_{23}V_{24}$,
$T_{22}=V_{21}V_{23}V_{24}V_{25}$,
$T_{23}=V_{20}V_{22}V_{23}V_{24}$,
$T_{24}=V_{22}V_{23}V_{24}V_{26}$,
$T_{25}=V_{23}V_{24}V_{25}V_{27}$, and
$T_{26}=V_{23}V_{24}V_{26}V_{27}$ (see on the figure at the
middle). All listed tetrahedra are taken by some integer affine
transformations to the unit basis tetrahedron. The integer
distance from the origin to the planes containing the faces
$T_{21},\ldots, T_{26}$ equal 1, 2, 2, 4, 8, and 13, respectively.
\end{example}

\begin{example} Consider the operator $A_3=A_{-1,-3,1,3}$. The group
$\Xi(A_3)$ is generated by the operators $B_{31}=A_3^{-2}$,
$B_{32}=(A_3{-}E)A_3^{-1}$, and $B_{33}=A_3{+}E$. Any fundamental
domain of the sail with $(0,0,0,1)$ as a vertex contains a unique
three-dimensional face. The polyhedron
$V_{30}V_{31}V_{32}V_{33}V_{34}V_{35}V_{36}V_{37}$ shown on the
figure (the right one) is an example of such face, here
$V_{30}=(-1, -1, -1, 0)$, $V_{31}=B_{33}(V_{30})$,
$V_{32}=B_{32}B_{33}(V_{30})$, $V_{33}=B_{31}B_{32}^{-1}(V_{30})$,
$V_{34}=B_{32}^{-1}(V_{30})$, $V_{35}=B_{31}B_{33}^{2}(V_{30})$,
$V_{36}=B_{31}B_{33}(V_{30})$,
$V_{37}=B_{31}B_{32}^{-1}B_{33}(V_{30})$.
The described face is contained in the plane on the unit distance
from the origin. The integer volume of the face equals $8$.
\end{example}

Example~3 provides the negative answer to the following question
for the case of $n=3$: {\it is it true, that any $n$-periodic
$n$-dimensional sail contains an $n$-dimensional face in some
hyperplane on integer distance to the origin greater than one?}
The answers for $n=2,4,5,6,\ldots$ are unknown. The answer to the
following question is also unknown to the author: {\it is it
true, that any $n$-periodic $n$-dimensional sail contains an
$n$-dimensional face in some hyperplane on unit integer distance
to the origin?}

$$
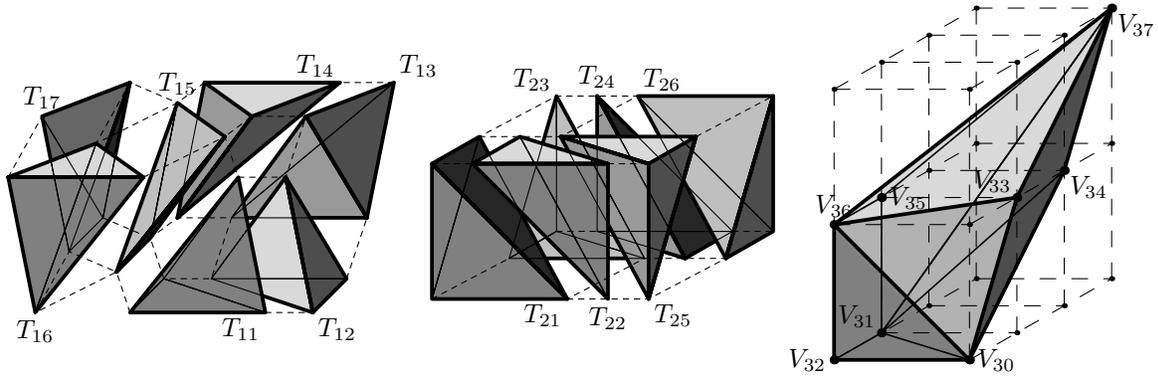
\epsfbox{examples.1}$$

We show with dotted lines how to glue the faces to obtain the
combinatorial scheme of the described fundamental domains.

\end{document}